\documentclass{elsart}

\usepackage{amssymb}
\usepackage{amsmath}

\newcommand{\lc}{l.c.}
\newcommand{\nsf}{n.s.f.}
\newcommand{\quotient}[2]{\raisebox{-0.195cm}[0.38cm]{}\mbox{\raisebox{-1pt}  
{$#1/\hspace{-2pt}\raisebox{-0.05cm}{${\textstyle#2}$}
$}}}
\def\11{\mbox1\hspace{-.25em}\text{I}}
\journal{Journal of Functional Analysis}

\begin{document}
\begin{frontmatter}
\title{On locally compact quantum groups\\
 whose algebras are factors}
\author{Pierre Fima}
\address{Laboratoire de math\'{e}matiques Nicolas Oresme, Universit\'{e} de Caen BP 5186, F14032 Caen Cedex, France}
\ead{pierre.fima@math.unicaen.fr}

\begin{abstract}
In this paper we are interested in examples of locally compact quantum groups $(M,\Delta)$ such that both von Neumann algebras, $M$ and the dual $\hat{M}$, are factors. There is a lot of known examples such that $(M,\hat{M})$ are respectively of type $(\rm{I}_{\infty},\rm{I}_{\infty})$ but there is no examples with factors of other types. We construct new examples of type $(\rm{I}_{\infty},\rm{II}_{\infty})$, $(\rm{II}_{\infty},\rm{II}_{\infty})$ and $(\rm{III}_{\lambda},\rm{III}_{\lambda})$ for each $\lambda\in [0,1]$. Also we show that there is no such example with $M$ or $\hat{M}$ a finite factor.
\end{abstract}
\begin{keyword}
Quantum groups \sep Factors \sep Crossed product
\MSC 46L52 \sep 46L65
\end{keyword}

\end{frontmatter}

\section{Introduction}
A locally compact (\lc) quantum group, in the von Neumann algebra setting (see \cite{KusVa2,KusVa3}), is a pair $(M,\Delta)$, where $M$ is a von Neumann algebra and $\Delta$ is a comultiplication on $M$, with left and right invariant weight. In a canonical way, every \lc{} group is a commutative \lc{} quantum group and every von Neumann group algebra of a \lc{} group is a cocommutative \lc{} quantum group. Conversely, every commutative or cocommutative \lc{} quantum group is obtained in this way. Our aim is to obtain new examples of \lc{} quantum groups which are as far as possible from groups so we will be interested in the "least commutative and cocommutative" examples. The formulation of the problem is the following. Given a pair of factors of a certain type $(x_{1}, x_{2})$, is it possible to find a \lc{} quantum group $(M,\Delta)$ such that $M$ is a type $x_{1}$ factor and the dual algebra $\hat{M}$ is a type $x_{2}$ factor (in the sense of Murray-von Neumann's and Connes
 ' classification of factors) ? There exists a lot of examples for the case $(\rm{I}_{\infty},\rm{I}_{\infty})$. In particular, it is shown in \cite{BaaSkVa} that for any matched pair of conjugated \lc{} groups the bicrossed product is a type $(\rm{I}_{\infty},\rm{I}_{\infty})$ \lc{} quantum group.

We start with a negative result showing that, if $(M,\Delta)$ is a \lc{} quantum group and $M$ is a finite factor, then $\hat{M}$ is not a factor. On a positive side, our tool to construct examples of \lc{} quantum groups is to use the bicrossed product of \lc{} groups (see \cite{VaesVain}). One can show that, if an action of a \lc{} group on a von Neumann algebra is free on its center, then the type of the cocycle crossed product does not depend a lot on the cocycle. This is why we only consider trivial cocycles. We construct for any $\lambda\in [0,1]$ a \lc{} quantum group $(M,\Delta)$ such that $M$ and $\hat{M}$ are type $\rm{III}_{\lambda}$ factors. Similary we obtain examples of type $(\rm{I}_{\infty},\rm{II}_{\infty})$ and $(\rm{II}_{\infty},\rm{II}_{\infty})$. All of them are ITPFI factors, a kind of infinite tensor product of the p-adic version of the Baaj and Skandalis' example (see \cite{VaesVain}).

This paper is organized as follows. In the second section we introduce some notations and recall some elementary facts about \lc{} quantum groups, bicrossed product construction, infinite tensor product and ITPFI factors. In the third section we prove that a \lc{} quantum group $(M,\Delta)$ with $M$ a finite factor is compact. In the fourth section we describe our examples using $p$-adic numbers.

\section{Preliminaries}
\paragraph{Locally compact quantum groups}

In this paper we suppose that all von Neumann algebras have separable predual and \lc{} groups are second countable. We denote by $\otimes$ the tensor product of Hilbert spaces or von Neumann algebras. We refer to \cite{Str2} for the theory of normal semifinite faithful (\nsf) weights on von Neumann algebras. If $\varphi$ is a \nsf{} weight on $M$, we use the standard notation
\begin{displaymath}
\mathcal{M}_{\varphi}^{+}=\{x\in M^{+}\,|\,\varphi(x)<\infty\},\,\,
\mathcal{N}_{\varphi}=\{x\in M\,|\,x^{*}x\in\mathcal{M}_{\varphi}^{+}\},
\,\, \mathcal{M}_{\varphi}=\mathcal{N}_{\varphi}^{*}\mathcal{N}_{\varphi}.
\end{displaymath}
We use \lc{} quantum groups in the von Neumann algebraic setting (see \cite{KusVa2}). A pair $(M,\Delta)$ is called a \lc{} quantum group when
\begin{itemize}
\item $M$ is a von Neumann algebra and $\Delta\rightarrow M\otimes M$ is a normal and unital *-homomorphism satisfying the coassociativity condition : 
\begin{displaymath}
(\Delta\otimes\iota)\Delta=(\iota\otimes\Delta)\Delta,
\end{displaymath}
where $\iota$ is the identity map.
\item There exist \nsf{} weights $\varphi$ and $\psi$ on $M$ such that
\begin{itemize}
\item $\varphi$ is left invariant in the sense that 
\begin{displaymath}
\varphi\left((\omega\otimes\iota)\Delta(x)\right)=\varphi(x)\omega(1),\,\,\forall x\in\mathcal{M}_{\varphi}^{+},\,\,\forall \omega\in M_{*}^{+},
\end{displaymath}
\item $\psi$ is right invariant in the sense that 
\begin{displaymath}
\psi\left((\iota\otimes\omega)\Delta(x)\right)=\psi(x)\omega(1),\,\,\forall x\in\mathcal{M}_{\psi}^{+},\,\,\forall\omega\in M_{*}^{+}.
\end{displaymath}
\end{itemize}
\end{itemize}
From \cite{KusVa2} we know that left invariant weights on $(M,\Delta)$ are unique to a positive scalar and the same holds for right invariant weights.

A \lc{} quantum group is called compact if its left invariant weight is finite.

Let $(M,\Delta)$ be a \lc{} quantum group, fix a left invariant \nsf{} weight $\varphi$ on $(M,\Delta)$ and represent $M$ on the GNS-space of $\varphi$ such that $(H,\iota,\Lambda)$ is a GNS construction for $\varphi$. Then we can define a unitary $W$ on $H\otimes H$ by
\begin{displaymath}
W^{*}\left(\Lambda(a)\otimes\Lambda(b)\right)=(\Lambda\otimes\Lambda)\left(\Delta(b)(a\otimes 1)\right)\qquad\text{for all  }a,b\in\mathcal{N}_{\varphi}.
\end{displaymath}
Where $\Lambda\otimes\Lambda$ is the canonical GNS-map for the tensor product weight $\varphi\otimes\varphi$. $W$ is called the fondamental unitary of $(M,\Delta)$. The comultiplication can be given in terms of $W$ by the formula $\Delta(x)=W^{*}(1\otimes x)W$ for all $x\in M$. Also the von Neumann algebra $M$ can be written in terms of $W$ as
\begin{displaymath}
M=\{(\iota\otimes\omega)(W)|\omega\in\mathcal{B}(H)_{*}\}^{-\sigma\mathrm{-strong}*},
\end{displaymath}
where $\{X\}^{-\sigma\mathrm{-strong}*}$ denote the $\sigma$-strong* closure of $X$. It is possible to define a new von Neumann algebra
\begin{displaymath}
\hat{M}=\{(\omega\otimes\iota)(W)|\omega\in\mathcal{B}(H)_{*}\}^{-\sigma\mathrm{-strong}*},
\end{displaymath}
and a comultiplication on $\hat{M}$ by $\hat{\Delta}(x)=\Sigma W(x\otimes 1)W^{*}\Sigma$ for $x\in\hat{M}$, where $\Sigma$ is the flip map on $H\otimes H$. Also, one can construct left and right invariant weight on $\hat{M}$ for $\hat{\Delta}$. We obtain in this way a new \lc{} quantum group $(\hat{M},\hat{\Delta})$ called the dual of $(M,\Delta)$. From \cite{KusVa2} we know that the bidual quantum group $(\hat{\hat{M}},\hat{\hat{\Delta}})$ is isomorphic to $(M,\Delta)$.

If $G$ is an ordinary locally compact group then $M=L^{\infty}(G)$ with the comultiplication $\Delta_{G}(f)(s,t)=f(st)$ and left and right invariant weight obtained by left and right Haar measure is a commutative \lc{} quantum group. Every commutative \lc{} quantum group is obtained in this way. Now, take $\hat{M}=\mathcal{L}(G)$ the group von Neumann algebra and $\hat{\Delta}_{G}(\lambda_{g})=\lambda_{g}\otimes\lambda_{g}$ where $(\lambda_{g})$ is the left regular representation of $G$. One can show that the Plancherel weight on $\hat{M}$ is left invariant. Also, it is easy to obtain a right invariant weight on $\hat{M}$ so $(\mathcal{L}(G),\hat{\Delta}_{G})$ is a \lc{} quantum group, this the dual of $(L^{\infty}(G),\Delta_{G})$. It is obvious that $\sigma\hat{\Delta}_{G}=\hat{\Delta}_{G}$ where $\sigma$ is the flip map on $M\otimes M$.  We say that $(\mathcal{L}(G),\hat{\Delta}_{G})$ is cocommutative. One can show that every cocommutative \lc{} quantum group is obtained in this way.

\paragraph{Bicrossed product}

Let $G_{1}, G_{2}$ be two closed subgroups of a \lc{} group $G$ such that $G_{1}\cap G_{2}=\{e\}$ and  $\mu(G-G_{1}G_{2})=0$ where $e$ is the identity element of $G$ and $\mu$ is a Haar measure on $G$. We say that the pair $(G_{1}, G_{2})$ is matched. We will now describe the bicrossed product construction of \lc{} group (see \cite{VaesVain}). Let $(G_{1}, G_{2})$ be a matched pair of \lc{} groups and $g\in G_{1}$, $s\in G_{2}$, then we can write nearly everywhere
\begin{displaymath}
gs=\alpha_{g}(s)\beta_{s}(g).
\end{displaymath}
We obtain two maps defined nearly everywhere and measurable
\begin{align*}
\alpha\,:\,G_{1}\times G_{2}\rightarrow G_{2}\,:\qquad (g,s)\rightarrow\alpha_{g}(s), \\
\beta\,:\,G_{2}\times G_{1}\rightarrow G_{1}\,:\qquad (s,g)\rightarrow\beta_{s}(g).
\end{align*}
Now we define two normal unital *-homomorphisms 
\begin{align*}
\alpha\, : \, L^{\infty}(G_{2})\rightarrow L^{\infty}(G_{1}\times G_{2})\,\, : \,\, (\alpha f)(g,s)=f(\alpha_{g}(s)),\\
\beta\, : \, L^{\infty}(G_{1})\rightarrow L^{\infty}(G_{2}\times G_{1})\,\, : \,\, (\beta f)(s,g)=f(\beta_{s}(g)),
\end{align*}
and we have 
\begin{displaymath}
(\iota\otimes\alpha)\alpha=(\Delta_{G_{1}}\otimes\iota)\alpha\qquad\text{and}\qquad
(\iota\otimes\beta)\beta=(\Delta_{G_{2}}\otimes\iota)\beta.
\end{displaymath}
Hence $\alpha$ will be an action of $G_{1}$ on the von Neumann algebra $L^{\infty}(G_{2})$ and $\beta$ an action of $G_{2}$ on the von Neumann algebra $L^{\infty}(G_{1})$. So we can define the crossed product von Neumann algebra $M=G_{1}\ltimes L^{\infty}(G_{2})$ and a faithfull *-homomorphism $\Delta$ : $M\rightarrow M\otimes M$
\begin{displaymath}
\Delta(x)=W^{*}(1_{L^{2}(G_{1}\times G_{2})}\otimes x)W\qquad\text{for }x\in M,
\end{displaymath}
where $W$ is a unitary in $L^{2}(G_{1}\times G_{2}\times G_{1}\times G_{2})$ defined by 
\begin{displaymath}
(W\xi)(g,s,h,t)=\xi(\beta_{\alpha_{g}(s)^{-1}t}(h)g,s,h,\alpha_{g}(s)^{-1}t).
\end{displaymath}
Then one can prove that $\Delta$ is a comultiplication on $M$ and the dual weight of the left invariant integral on $L^{\infty}(G_{2})$ is left invariant for $\Delta$. In this way we obtain a \lc{} quantum group $(M,\Delta)$ with dual $(\hat{M},\hat{\Delta})$ such that 
\begin{displaymath}
\hat{M}=L^{\infty}(G_{1})\rtimes G_{2},\quad
\hat{\Delta}(x)=\Sigma W(x\otimes 1_{L^{2}(G_{1}\times G_{2})})W^{*}\Sigma\quad\text{for }x\in \hat{M}.
\end{displaymath}

\paragraph{Infinite tensor product of von Neumann algebras}

For each $n\in\Nset$ let $M_{n}$ be a von Neumann algebra acting on an Hilbert space $H_{n}$ and $\xi_{n}$ a norm $1$ vector in $H_{n}$. The infinite tensor product of $M_{n}$ relatively to $\xi_{n}$ is the von Neumann algebra generated by the operators $x_{1}\otimes\ldots\otimes x_{k}\otimes 1\otimes\ldots$ for $k\in \Nset$ and $x_{i}\in
 M_{i}$ in the infinite tensor product of Hilbert spaces $H_{n}$ relatively to the norm $1$ vectors $\xi_{n}$. We denote this von Neumann algebra by  
\begin{displaymath}
\bigotimes(M_{n},H_{n},\xi_{n}).
\end{displaymath}
If each $M_{n}$ is a factor then $\bigotimes(M_{n},H_{n},\xi_{n})$ is a factor (see \cite{ArW}).

The following lemma is certainly well known but we could not find a proof in the literature.
\begin{lem}\label{valprop}
Let $p_{n}\neq 0$ be a projection in a von Neumann algebra $M_{n}$ and $\omega_{n}$ a normal faithful state on $M_{n}$ with GNS space $H_{n}$ such that $\omega_{n}=\omega_{\xi_{n}}$. Put $M=\bigotimes\left(M_{n},\omega_{n}\right)$. The decreasing sequence of projections
\begin{displaymath}
p_{1}\otimes\ldots\otimes p_{n}\otimes 1\ldots
\end{displaymath}  
converge to a projection $p\in M$ and we have 
\begin{displaymath}
p\neq 0\Leftrightarrow\sum_{n}||(1-p_{n})\xi_{n}||^{2}<\infty,
\end{displaymath}
and if $p\neq 0$
\begin{displaymath}
pMp\simeq
\bigotimes\left(p_{n}M_{n}p_{n},\omega_{\eta_{n}}\right),\quad\text{where}\quad \eta_{n}:=\frac{p_{n}\xi_{n}}{||p_{n}\xi_{n}||}.
\end{displaymath}
\end{lem}
\begin{pf}
$p\neq 0$ if and only if $(\otimes\omega_{n})(p)>0$. This is equivalent to 
\begin{displaymath}
\sum_{n}||(1-p_{n})\xi_{n}||^{2}<\infty.
\end{displaymath}
The isomorphism is a simple identification.\qed
\end{pf}

\begin{exmp}
For each $n\in\Nset$ let $M_{n}$ be a von Neumann algebra, $\varphi_{n}$ a \nsf{} weight on $M_{n}$ and $q_{n}\in M_{n}$ a projector with $\varphi_{n}(q_{n})=1$. We can take $H_{n}$ the Hilbert space of the G.N.S. construction for $\varphi_{n}$ and $\xi_{n}=\Lambda_{\varphi_{n}}(q_{n})$. We introduce the notation 
\begin{displaymath}
M=\bigotimes\left(M_{n},\varphi_{n},q_{n}\right):= \bigotimes(M_{n},H_{n},\xi_{n}).
\end{displaymath}
Observe that there is a natural projector $q\in M$ defined by the infinite tensor product of the $q_{n}$. When $\varphi_{n}=\omega_{n}$ is a normal faithful state on $M_{n}$ and $q_{n}=1$, we use the standard notation 
\begin{displaymath}
\bigotimes\left(M_{n},\omega_{n}\right):=\bigotimes\left(M_{n},\omega_{n},1\right).
\end{displaymath}
\end{exmp}

\begin{prop}\label{qMq}
Let $M_{n}$, $\varphi_{n}$ and $q_{n}$ be as above and suppose that $q_{n}\in M_{n}^{\varphi_{n}}$. Denote by $\omega_{\xi_{n}}$ the vector state associated to $\xi_{n}$. Then $\omega_{\xi_{n}}$ is faithful on $q_{n}M_{n}q_{n}$ and
\begin{displaymath}
qMq\simeq\bigotimes\left(q_{n}M_{n}q_{n},\omega_{\xi_{n}}\right).
\end{displaymath}
\end{prop}

\begin{pf}
Let $J_{n}$ be the usual antiunitary operator associated with $\varphi_{n}$ and $J=\otimes J_{n}$, $q=\otimes q_{n}$. We have $||x\xi_{n}||^{2}=\varphi_{n}(q_{n}x^{*}xq_{n})$ thus the faithfulness of $\varphi_{n}$ on $M_{n}$ implies that $\omega_{\xi_{n}}$ is faithful on $q_{n}M_{n}q_{n}$. Note that the close linear subspace of $q_{n}H_{n}$ generated by $q_{n}M_{n}q_{n}\xi_{n}$ is $J_{n}q_{n}J_{n}q_{n}H_{n}$. Thus, the GNS space of $\omega:=\otimes\omega_{\xi_{n}}$ on $N:=\bigotimes\left(q_{n}M_{n}q_{n},\omega_{\xi_{n}}\right)$ is canonically isomorphic with $qJqJH$, and the image of $qMq$ by the restriction homomorphism to the invariant subspace $qJqJ H$ is $N$. This homomorphism is in fact an isomorphism because the closure of $M^{'}qJqJH$ is $qH$. \qed
\end{pf}

\begin{rem}
One can show that if $q_{n}\in M_{n}^{\varphi_{n}}$ there exists a canonical \nsf{} weight $\varphi$ on $M$ such that $\sigma_{t}^{\varphi}=\otimes_{n}\sigma_{t}^{\varphi_{n}}$ and $\varphi q =\otimes_{n} \varphi_{n} q_{n}$. This is the noncommutative analogue of the restricted direct product of measurable spaces with non necessarily finite measure.
\end{rem} 

\paragraph{Resctricted direct product action}

In the sequel, all group actions on von Neumann algebras are
supposed to be pointwise $\sigma$-weakly continuous. Let $G_{n}$ be
a sequence of \lc{} groups and $\mu_{n}$ a left Haar measure on
$G_{n}$. We suppose that for all $n$ there is a compact open
subgroup $K_{n}$ of $G_{n}$ such that $\mu_{n}(K_{n})=1$. Recall
that the restricted direct product $\prod\,^{'}(G_{n},K_{n})$ is
defined as the set of $(x_{n})\in\prod G_{n}$ such that $x_{n}\in
K_{n}$ for $n$ large enough (see \cite{Bl} for details). Let $M_{n}$
be a sequence of von Neumann algebras with actions 
$\alpha^{n}\,:\,G_{n}\rightarrow \text{Aut}(M_{n})$. Let
$\varphi_{n}$ be a \nsf{} weight on $M_{n}$ and $q_{n}$ a projection
in the centralizer of $\varphi_{n}$ with $\varphi_{n}(q_{n})=1$ and such that for all $n$, for all
$g_{n}\in G_{n}$, there exists $c_{n}(g_{n})>0$ such that
$\varphi_{n}\circ\alpha^{n}=c_{n}(g_{n})\varphi_{n}$,
$c_{n}|_{K_{n}}=1$ and for all $g_{n}\in K_{n}$ one has
$\alpha^{n}_{g_{n}}(q_{n})=q_{n}$. With this data one can construct
an action of the restricted direct product
$G=\prod^{'}\left(G_{n},K_{n}\right)$ on the infinite tensor product $M=\bigotimes\left(M_{n},\varphi_{n},q_{n}\right)$. We fix a G.N.S. construction $(H_{n},\iota,\Lambda_{n})$ for $\varphi_{n}$ and we put $H=\bigotimes\left(H_{n},\Lambda_{n}(q_{n})\right)$.

\begin{prop}
There exists a unique action $\alpha\,:\,G\rightarrow\text{Aut}(M)$, called the restricted direct product of $\alpha^{n}$, such that for all $g=(g_{n})\in G$ and $x_{n}\in M_{n}$
\begin{displaymath}
\alpha_{g}(x_{1}\otimes\ldots\otimes x_{n}\otimes 1\ldots)
=\alpha^{1}_{g_{1}}(x_{1})\otimes\ldots\otimes\alpha^{n}_{g_{n}}(x_{n})\otimes 1\otimes\ldots
\end{displaymath}
\end{prop}

\begin{pf}
The uniqueness part is obvious. To show the existence, we first compute a unitary implementation of the actions $\alpha^{n}$. It is easy to see that, for $g_{n}\in G_{n}$, the operator
\begin{displaymath}
U_{n}(g_{n})\,:\,\Lambda_{n}(x)\mapsto c_{n}(g_{n})^{\frac{1}{2}}\Lambda_{n}(\alpha^{n}_{g_{n}}(x)),\qquad\text{for }x\in\mathcal{N}_{\varphi_{n}}, 
\end{displaymath}
can be extended to a unitary operator on $H_{n}$, still denoted by $U_{n}(g_{n})$, and such that
\begin{displaymath}
\alpha^{n}_{g_{n}}(x)=U_{n}(g_{n})xU_{n}(g_{n})^{*}\qquad\text{for all }x\in M_{n}.
\end{displaymath}
The hypotesis implies that for all $n$ and for all $g_{n}\in K_{n}$ one has $U_{n}(g_{n})\Lambda_{n}(q_{n})=\Lambda_{n}(q_{n})$. Then for all $g\in G$ one can define a unitary operator $U_{g}$ on $H$ by $U_{g}=\otimes_{n} U_{n}(g_{n})$ where $g=(g_{n})$. In this way we obtain a group homomorphism $g\mapsto U_{g}$ from $G$ to the unitary group of $H$. Because we have $U_{g} M U_{g}^{*}=M$ this allows us to construct a group homomorphism $\alpha\,:G\rightarrow\text{Aut}(M)$ defined by
\begin{displaymath}
\alpha_{g}(x)=U_{g} x U_{g}^{*}\qquad\text{for }x\in M.
\end{displaymath}
This is obvious that $\alpha$ is pointwise $\sigma$-weakly continuous and verifies the equation.\qed
\end{pf}

Let us identify the crossed product of $G$ by $M$ with an infinite tensor product of the crossed products of $G_{n}$ by $M_{n}$. We denote by $\pi_{n}$ the inclusion of $M_{n}$ into $G_{n}\ltimes M_{n}$ and $\pi$ the inclusion of $M$ into $G\ltimes M$. We denote by $\tilde{\varphi}_{n}$ the dual weight of $\varphi_{n}$ and by $\11_{A}$ the caracteristic function of a mesurable set $A$.

\begin{prop}\label{crossed}
Let $e_{n}=\pi_{n}(q_{n})\left(\lambda(\11_{K_{n}})\otimes 1\right)$, where $\lambda(\11_{K_{n}})$ is the convolution operator by $\11_{K_{n}}$. Then $e_{n}$ is a projection in $G_{n}\ltimes M_{n}$. Moreover, one has $\tilde{\varphi}_{n}(e_{n})=1$, $e_{n}\in \left(G_{n}\ltimes M_{n}\right)^{\tilde{\varphi}_{n}}$ and
\begin{displaymath}
G\ltimes M\simeq
\bigotimes\left(G_{n}\ltimes M_{n},\tilde{\varphi}_{n},e_{n}\right).
\end{displaymath}
\end{prop}

\begin{pf}
Because $\pi_{n}(q_{n})$ and $\lambda(\11_{K_{n}})$ are projections, if $\pi_{n}(q_{n})$ and $\lambda(\11_{K_{n}})\otimes 1$ commute then $e_{n}$ is a projection. Take $\xi\in L^{2}(G_{n},H_{n})$ then
\begin{eqnarray*}
\lefteqn{\left((\lambda(\11_{K_{n}})\otimes 1)\pi_{n}(q_{n})\xi\right)(g)=\int_{G_{n}}\11_{K_{n}}(t)\alpha^{n}_{g^{-1}t}(q_{n})\xi(t^{-1}g)\d\mu_{n}(t)}\\
& & =\alpha^{n}_{g^{-1}}(q_{n})\int_{G_{n}}\11_{K_{n}}(t)\xi(t^{-1}g)\d\mu_{n}(t),\,\,\,\text{because}\,\,\,\forall t\in K_{n},\,\,\,\alpha^{n}_{t}(q_{n})=q_{n}\\
& & =\left(\pi_{n}(q_{n})(\lambda(\11_{K_{n}})\otimes1)\xi\right)(g).
\end{eqnarray*}
Thus $e_{n}$ is a projection. Now, using $K_{n}\subset\text{Ker}(\delta_{G_{n}})$ and $c_{n}|_{K_{n}}=1$, we have
\begin{displaymath}
\sigma_{t}^{\tilde{\varphi}_{n}}(\lambda(\11_{K_{n}})\otimes1)=\lambda(\11_{K_{n}})\otimes1.
\end{displaymath}
This implies that
\begin{displaymath}
\sigma_{t}^{\tilde{\varphi}_{n}}(e_{n})=\pi_{n}(\sigma_{t}^{\varphi_{n}}(q_{n}))
\sigma_{t}^{\tilde{\varphi}_{n}}(\lambda(\11_{K_{n}})\otimes1)=e_{n}.
\end{displaymath}
Next, using definition of the dual weight, we have $\tilde{\varphi}_{n}(e_{n})=\varphi_{n}(q_{n})\11_{K_{n}}(1)
=1$. Recall that, using the classical explicit G.N.S. construction $(L^{2}(G_{n},M_{n}),\iota,\tilde{\Lambda}_{n})$ for the dual weight, one has (see \cite{Vaes})
\begin{displaymath}
\tilde{\Lambda}_{n}(e_{n})=\11_{K_{n}}\otimes\Lambda_{n}(q_{n}).
\end{displaymath}
We denote this vector by $\xi_{n}$. We define the operator
\begin{displaymath}
U\,:\,\bigotimes\left(L^{2}(G_{n},H_{n}),\xi_{n}\right)
\rightarrow L^{2}(G,H)
\end{displaymath}
on a dense subset by
\begin{displaymath}
U(F_{1}\otimes\ldots\otimes F_{n}\otimes\bar{\xi}_{n})(g)
=\bigotimes_{i=1}^{n}F_{i}(g_{i})\otimes\left(
\bigotimes_{i=n+1}^{\infty}\11_{K_{i}}(g_{i})\Lambda_{i}(p_{i})\right),
\end{displaymath}
where $F_{i}\in L^{2}(G_{i},H_{i})$, $g=(g_{n})\in G$ and $\bar{\xi}_{n}=\otimes_{i=n+1}^{\infty}\xi_{i}$. Then $U$ is an isometry with dense range. Thus we obtain a unitary operator, again denoted by $U$, such that, if $g=(g_{1},\ldots,g_{n},1,\ldots)$,
\begin{displaymath}
U\left((\lambda_{g_{1}}\otimes1)\otimes\ldots\otimes(\lambda_{g_{n}}\otimes 1)\otimes 1\otimes\ldots\right) U^{*}
=\lambda_{g}\otimes 1,
\end{displaymath}
\begin{displaymath}
U\left(\pi_{1}(x_{1})\otimes\ldots\otimes\pi_{n}(x_{n})\otimes 1\otimes\ldots\right) U^{*}
=\pi(x_{1}\otimes\ldots\otimes x_{n}\otimes 1\otimes\ldots).
\end{displaymath}
It follows that
\begin{displaymath}
U\left(\bigotimes\left(G_{n}\ltimes M_{n},\tilde{\varphi}_{n},e_{n}\right)\right) U^{*}=G\ltimes M.\qed
\end{displaymath}
\end{pf}

\paragraph{ITPFI factors and Boca-Zaharescu factors}

In \cite{ArW} Araki and Woods define ITPFI factors as infinite tensor product of type $\rm{I}$ factors
\begin{displaymath}
\bigotimes(M_{n},H_{n},\xi_{n}),
\end{displaymath}
where $M_{n}$ is a type $\rm{I}$ factor acting on $H_{n}$ and $\xi_{n}$ is a norm $1$ vector in $H_{n}$. All these factors are hyperfinite.

If $M$ is a type $\rm{I}$ factor acting on $H$, we can write $H=H_{1}\otimes H_{2}$ such that $M=\mathcal{B}(H_{1})\otimes 1$. Now, let $\Omega\in H$ be a norm $1$ vector and consider the normal state on $M$
\begin{displaymath}
\omega(x)=\langle (x\otimes 1)\Omega,\Omega \rangle.
\end{displaymath}
Hence there exists a density matrix  $\rho_{\Omega}\in\mathcal{B}(H_{1})$ such that $\omega(x)=Tr(\rho_{\Omega}x)$. It is easy to see that the ordered list  (with multiplicity) of the non zero eigenvalues of the operator $\rho_{\Omega}$ does not depend on the decomposition of $H$ in $H_{1}\otimes H_{2}$. This list is denoted by Sp$(\Omega\,|\,M)$. The type of the ITPFI factor $\bigotimes\left(M_{n},H_{n},\xi_{n}\right)$ only depends on the list Sp$(\xi_{n}\,|\,M_{n})$. In the fourth section we will use the fact that if each $M_{n}$ is a type $\rm{I}_{n_{\nu}}$ factor, with $2\leq n_{\nu}\leq\infty$, and Sp$(\xi_{n}\,|\,M_{n})=\{\lambda_{n_{i}},\,i=1,2,\ldots,n_{\nu}\}$ then, if $\lambda_{n_{1}}\geq\delta$ for some $\delta>0$ and for all $n$, $M$ is a type $\rm{III}$ factor if and only if
\begin{displaymath}
\sum_{n,i}\lambda_{n_{i}}\inf\left\{\left|\frac{\lambda_{n_{1}}}{\lambda_{n_{i}}}-1\right|^{2},C\right\}=\infty,
\end{displaymath}
for some positive $C$.

Let $\mathcal{S}$ be an infinite subset of the set $\mathcal{P}$ of all prime numbers and $\beta\in ]0,1]$. In \cite{Boca} Boca and Zaharescu studied the following ITPFI factor
\begin{displaymath}
M_{\beta,\mathcal{S}}:=\bigotimes_{p\in\mathcal{S}}\left(\mathcal{B}(l^{2}(\Nset)),\omega_{p,\beta}\right),
\end{displaymath}
where $\omega_{p}(x):=\sum_{n}p^{-n\beta}(1-p^{-\beta})\langle xe_{n},e_{n}\rangle$ and $(e_{n})$ is the canonical basis of $l^{2}(\Nset)$. We denote by $N_{\mathcal{S}}$ the factor $M_{1,\mathcal{S}}$. In \cite{Boca} Boca and Zaharescu show that
\begin{enumerate}
\item  For any $\lambda\in [0,1]$ and $\beta\in ]0,1]$, there is a subset $\mathcal{S}$ of $\mathcal{P}$ such that $M_{\mathcal{\beta,S}}$ is a type $\rm{III}_{\lambda}$ factor.
\item  For any $\beta\in ]0,1]$, any countable subgroup $K$ of $\Rset$ and any countable subset $\Sigma$ of $\Rset -K$, there exists a subset $\mathcal{S}$ of $\mathcal{P}$ such that $T(M_{\beta,\mathcal{S}})$ contains $K$ and does not intersect $\Sigma$, 
\end{enumerate}
where $T(M)$ denotes the Connes' $T$ invariant of the von Neumann algebra $M$ (see \cite{Co1}).

\begin{rem}
It was shown in \cite{Boca} that $M_{\beta,\mathcal{S}}$ is an $\text{ITPFI}_{2}$ (infinite tensor product of type $\rm{I}_{2}$ factors) for all $\beta\in ]\frac{1}{2},1]$. In fact it is possible to show that for all $\beta\in ]0,1]$, $M_{\beta,\mathcal{S}}$ is an $\text{ITPFI}_{m}$ with $\beta>\frac{1}{m}$. Indeed, for such $m$ and $\beta$ put 
\begin{displaymath}
q_{p}(e_{n})=\left\{\begin{array}{cl}
e_{n} & \text{if $0\leq n \leq m-1$}\\
0 & \text{otherwise.}\end{array}\right.
\end{displaymath} 
Then because of 
\begin{displaymath}
\sum_{p\in\mathcal{S}}\sum_{n\geq m}(1-p^{-\beta})p^{-n\beta}=\sum_{p\in\mathcal{S}}p^{-m\beta}<\infty,
\end{displaymath}
for all $\beta>\frac{1}{m}$, we can apply Lemma \ref{valprop} to obtain a projection $p\neq 0$ such that $p\left(M_{\beta,\mathcal{S}}\right)p$ is an $\text{ITPFI}_{m}$. Moreover, it is easy to see that $p$ is purely infinite, thus $M_{\beta,\mathcal{S}}$ is $\text{ITPFI}_{m}$.
\end{rem}
\section{The case of a finite factor}

In this section we show that if $(M,\Delta)$ is a \lc{} quantum group such that $M$ is a finite factor, then $(M,\Delta)$ is compact so $\hat{M}$, being an infinite direct sum of full matrix algebras, is not a factor.

The idea of the proof of the next lemma was taken from \cite{vanDaele}.
\begin{lem}
Let $(M,\Delta)$ be a \lc{} quantum group. Suppose that $M$ is a finite factor. Let $\tau$ be the unique tracial state on $M$. Then, for all $\rho\in M_{*}$ with $0\leq \rho\leq\tau$, one has :
\begin{displaymath}
\rho * \tau=\tau * \rho=\rho(1)\tau.
\end{displaymath}
\end{lem}

\begin{pf}
Let $a$ be in $M$ and define $b=(\iota\otimes\tau)\Delta(a)$. Then, by unicity of $\tau$, one has $\tau * \tau = \tau$, and using the coassociativity of $\Delta$ we obtain
\begin{align*}
(\iota\otimes \tau)\Delta(b)
 & =(\iota\otimes\tau)\Delta((\iota\otimes\tau)\Delta(a))
=(\iota\otimes\tau\otimes\tau)((\Delta\otimes\iota)\Delta(a)) \\
 & =(\iota\otimes\tau\otimes\tau)((\iota\otimes\Delta)\Delta(a))\\
 & =(\iota\otimes(\tau * \tau))\Delta(a)=(\iota\otimes\tau)\Delta(a)=b.
\end{align*}

This implies the following relations.
\begin{align}
(\iota\otimes\tau)\left((b^{*}\otimes 1)\Delta(b)\right)
=b^{*}(\iota\otimes\tau)\Delta(b)=b^{*}b\label{eqtypII11}, \\
(\iota\otimes\tau)\left(\Delta(b^{*})(b\otimes 1)\right)
=\left((\iota\otimes 1\Delta(b)\right)^{*} b=b^{*}b.\label{eqtypeII12}
\end{align}

Now define
\begin{align*}
k
 & =\left(\Delta(b)-b\otimes 1\right)^{*}\left(\Delta(b)-b\otimes 1\right) \\
 & =\Delta(b^{*}b)-(b^{*}\otimes 1)\Delta(b)-\Delta(b^{*})(b\otimes 1)+b^{*}b\otimes 1.
\end{align*}

Then $k\geq 0$ and, from the equations (\ref{eqtypII11}) and (\ref{eqtypeII12}), we obtain
\begin{align*}
(\tau\otimes\tau)(k)
 & =\tau((\iota\otimes\tau)(k))=(\tau\otimes\tau)(\Delta(b^{*}b)-\tau(b^{*}b) \\
 & =\tau * \tau(b^{*}b)-\tau(b^{*}b)=0.
\end{align*}
Then, if $\rho\in M_{*}$ with $0\leq \rho\leq\tau$, one has $(\tau\otimes\rho)(k)\leq (\tau\otimes\tau)(k)=0$. This implies, with the Cauchy-Schwartz inequality, that for all $c\in M$ we have
\begin{displaymath}
(\tau\otimes\rho)\left((c\otimes 1)(\Delta(b)-b\otimes 1)\right)=0\qquad\text{thus},
\end{displaymath}
\begin{displaymath}
(\tau\otimes\rho)\left((c\otimes 1)(\Delta(b)\right)=\rho(1)\tau(cb).
\end{displaymath}
Using the definition of $b$, we see that the last equation is equivalent to
\begin{align*}
 & (\tau\otimes\rho)\left((c\otimes 1)\Delta\left((\iota\otimes\tau)\Delta(a)\right)\right)=\rho(1)\tau\left(c(\iota\otimes\tau)\Delta(a)\right)\\
\Leftrightarrow
 & (\tau\otimes\rho\otimes\tau)\left((c\otimes 1\otimes 1)(\Delta\otimes\iota)\Delta(a)\right)=\rho(1)(\tau\otimes\tau)\left((c\otimes 1)\Delta(a)\right)\\
\Leftrightarrow
 & (\tau\otimes(\rho * \tau))\left((c\otimes 1)\Delta(a)\right)=\rho(1)(\tau\otimes\tau)\left((c\otimes 1)\Delta(a)\right),
\end{align*}
and this is true for all $a$ and $b$ in $M$. Now, because $\Delta(M)(M\otimes 1)$ is $\sigma$-weakly dense in $M\otimes M$ and $\tau$ is a trace we have, for all $x\in M\otimes M$, 
\begin{displaymath}
(\rho(1)\tau\otimes\tau)(x)=(\tau\otimes(\rho * \tau))(x).
\end{displaymath}

Putting $x=1\otimes y$ in the last equation, we obtain $\rho(1)\tau=\rho * \tau$. The proof of $\rho(1)\tau=\tau * \rho$ is the same.\qed
\end{pf}

We are now able to prove that a \lc{} quantum group $(M,\Delta)$ with $M$ a finite factor is compact.

\begin{thm}
Let $(M,\Delta)$ be a \lc{} quantum group with $M$ a finite factor. Then $(M,\Delta)$ is compact and $\tau$ is the Haar state on $M$, where $\tau$ is the unique tracial state on $M$.
\end{thm}

\begin{pf}
Let $(H,\Lambda,\iota)$ be a G.N.S. construction for $\tau$ and $J$ the canonical involutive isometry associated to $\tau$. Let $a$ be in $M$ and consider the positive normal linear form $\omega_{\Lambda(a)}$. We have
\begin{align*}
\omega_{\Lambda(a)}(x^{*}x)
 & =||\Lambda(xa)||^{2}=||Ja^{*}J\Lambda(x)||^{2} \\
 & \leq ||a||^{2}\tau(x^{*}x).
\end{align*}
This implies that $\frac{\omega_{\Lambda(a)}}{||a||^{2}}\leq\tau$ and, using the previous lemma, we conclude that $\omega_{\Lambda(a)} * \tau=\tau * \omega_{\Lambda(a)}=\omega_{\Lambda(a)}(1)\tau$ for all $a\in M$. Now, using that $M\subset\mathcal{B}(H)$ is standard, we know that if $\omega\in M_{*}$ and $\omega\geq 0$ there exists $\xi\in H$ such that $\omega=\omega_{\xi}$. Take a net $(a_{i})$ in $M$ such that $\Lambda(a_{i})$ converges in $H$ to $\xi$ then, for all $x\in M$, $\omega_{\Lambda(a_{i})}(x)$ converges to $\omega(x)$. In particular, for $x$ in $M$, we have
\begin{align*}
(\omega_{\Lambda(a_{i})} * \tau)(x)\rightarrow (\omega * \tau)(x) \\
(\tau * \omega_{\Lambda(a_{i})})(x)\rightarrow(\tau * \omega)(x).
\end{align*}
Because of
\begin{displaymath}
(\omega_{\Lambda(a_{i})} * \tau)(x)=(\tau * \omega_{\Lambda(a_{i})})(x)=||\Lambda(a_{i})||^{2}\tau(x)\rightarrow ||\xi||^{2}\tau(x)=\omega(1)\tau(x),
\end{displaymath}
we see that $\omega * \tau=\tau * \omega=\omega(1)\tau$ and, by linearity, the last equality holds for all $\omega\in M_{*}$. This concludes the proof.\qed
\end{pf}

\section{Examples}

Let $\mathcal{P}$ be the set of all prime numbers. In the sequel, if
$p$ is a prime number, we denote by $\Qset_{p}$ the field of rational
$p$-adic numbers and $\Zset_{p}$ the ring of $p$-adic integers. Let
$\mathcal{S}$ be an infinite subset of $\mathcal{P}$ and
$\mathcal{A}_{\mathcal{S}}$ the restricted direct product of $\Qset_{p}$ relatively to the compact open subgroups $\Zset_{p}$ for
$p\in\mathcal{S}$ (see \cite{Bl}) :
\begin{displaymath}
\mathcal{A}_{\mathcal{S}}=\prod_{p\in\mathcal{S}}\,^{'}\left(\Qset_{p},\Zset_{p}\right).
\end{displaymath}
Then $\mathcal{A}_{\mathcal{S}}$ is a second countable \lc{} ring.
The group of invertible elements of $\mathcal{A}_{\mathcal{S}}$ is
\begin{displaymath}
\mathcal{A}^{*}_{\mathcal{S}}=\prod_{p\in\mathcal{S}}\,^{'}\left(\Qset_{p}^{*},\Zset_{p}^{*}\right).
\end{displaymath}
Now, denote by $G_{\mathcal{S}}$ the $ax+b$-group of
$\mathcal{A}_{\mathcal{S}}$ :
\begin{displaymath}
G_{\mathcal{S}}=\mathcal{A}^{*}_{\mathcal{S}}\ltimes\mathcal{A}_{\mathcal{S}},
\end{displaymath}
and define the following subgroups.

\begin{displaymath}
G_{\mathcal{S}}^{1}=\left\{(a,0)\in G_{\mathcal{S}}\right\},\,\,\,
G_{\mathcal{S}}^{2}=\left\{\left((a_{p}),(b_{p})\right)\in
G,\,\,a_{p}+b_{p}p=1\,\,\forall p\in\mathcal{S}\right\}.
\end{displaymath}

We can rewrite $G_{\mathcal{S}}^{2}$ as
\begin{displaymath}
G_{\mathcal{S}}^{2}=\left\{
\begin{array}{c}
\left((a_{p}),\left(\frac{1-a_{p}}{p}\right)\right),\,\,a_{p}\neq 0\,\,\,\forall p\in\mathcal{S}\\
\text{  and  }a_{p}\in1+p\Zset_{p}\text{  for $p$ large enough}
\end{array}
\right\}.
\end{displaymath}
$G_{\mathcal{S}}^{1}$ is the subgroup of $G_{\mathcal{S}}$ which
fixes $0$. $G_{\mathcal{S}}^{2}$ is, formally, the subgroup of
$G_{\mathcal{S}}$ which fixes
$\left(\frac{1}{p}\right)_{p\in\mathcal{S}}$. We denote by
$\mu_{p}^{+}$ the additive Haar measure on $\Qset_{p}$ such that
$\mu_{p}^{+}(\Zset_{p})=1$ and by $\mu_{p}^{\times}$ the multiplicative
Haar measure on $\Qset_{p}^{*}$ such that $\mu_{p}^{\times}(\Zset_{p}^{*})=1$. Let
$\mu^{+}$ be the product measure of $\mu_{p}^{+}$, this is an
additive Haar measure on $\mathcal{A}_{\mathcal{S}}$, let
$\mu^{\times}$ be the product measure of $\mu_{p}^{\times}$, this is
a Haar measure on $\mathcal{A}^{*}_{\mathcal{S}}$. On
$G_{\mathcal{S}}$, the right Haar measure which is equal to $1$ on
$\prod_{p\in\mathcal{S}}\Zset_{p}^{*}\times\Zset_{p}$ is
$d\mu^{\times}(x)d\mu^{+}(y)$ and the left Haar measure which is
equal to $1$ on $\prod_{p\in\mathcal{S}}\Zset_{p}^{*}\times\Zset_{p}$ is
$\delta(x)d\mu^{\times}(x)d\mu^{+}(y)$, where
\begin{displaymath}
\delta(x)=\prod_{p\in\mathcal{S}} \frac{1}{|x_{p}|}_{p},\,\,\,x=(x_{p})\in\mathcal{A}_{\mathcal{S}}^{*}.
\end{displaymath}

We now prove the following easy lemma.

\begin{lem}
The groups $G_{\mathcal{S}}^{1}$, $G_{\mathcal{S}}^{2}$ are matched.
Moreover, the bicrossed product of $G_{\mathcal{S}}^{1}$ by
$G_{\mathcal{S}}^{2}$ is not regular, it is semi-regular in the
sense of \cite{BaaSkVa}.
\end{lem}

\begin{pf}
It is clear that $G_{\mathcal{S}}^{1}$ and $G_{\mathcal{S}}^{2}$ are
closed subgroups of $G_{\mathcal{S}}$ and $G_{\mathcal{S}}^{1}\cap
G_{\mathcal{S}}^{2}=\{1\}$. So we must prove that
$G_{\mathcal{S}}-G_{\mathcal{S}}^{1}G_{\mathcal{S}}^{2}$ is closed
and its Haar measure is zero. From
\begin{displaymath}
G_{\mathcal{S}}^{2}G_{\mathcal{S}}^{1}=
\left\{
\begin{array}{c}
\left((a_{p}b_{p}),\left(\frac{1-a_{p}}{p}\right)\right),
\,\,b=(b_{p})_{p\in\mathcal{S}}\in\mathcal{A}^{*}_{\mathcal{S}},\,\,a_{p}\neq 0\,\,\,\forall p\in\mathcal{S}\\
\text{  and  }a_{p}\in1+p\Zset_{p}\text{  for $p$ large enough}
\end{array}
\right\},
\end{displaymath}
we conclude that
\begin{displaymath}
G_{\mathcal{S}}^{2}G_{\mathcal{S}}^{1}=
\left\{(a,b)\in G_{\mathcal{S}},\,\,b_{p}\neq\frac{1}{p}\,\forall
p\in\mathcal{S}\text{  with  }
b=(b_{p})_{p\in\mathcal{S}}\right\}.
\end{displaymath}
It follows that
$G_{\mathcal{S}}^{2}G_{\mathcal{S}}^{1}$ is open and
\begin{displaymath}
G_{\mathcal{S}}-G_{\mathcal{S}}^{2}G_{\mathcal{S}}^{1}=\left\{(a,b)\in G_{\mathcal{S}},\,\,\exists p\in \mathcal{S},\,\,
b_{p}=\frac{1}{p}\text{  with }b=(b_{p})_{p\in\mathcal{S}}\right\}
\end{displaymath}
has Haar measure equal to zero.\qed
\end{pf}

Denote by $(M_{\mathcal{S}},\Delta_{\mathcal{S}})$ the bicrossed
product of $G_{\mathcal{S}}^{1}$ and $G_{\mathcal{S}}^{2}$.  Under the canonical identification of $G_{\mathcal{S}}^{1}$ with
$\mathcal{A}^{*}_{\mathcal{S}}$ and $G_{\mathcal{S}}^{2}$ with
$\mathcal{K}_{\mathcal{S}}$, where $\mathcal{K}_{\mathcal{S}}$ is
the following restricted direct product
\begin{displaymath}
\mathcal{K}_{\mathcal{S}}=\prod_{p\in\mathcal{S}}\,  ^{'}\left(\Qset_{p}^{*},1+p\Zset_{p}\right),
\end{displaymath}
the group actions $\alpha$ of $G_{\mathcal{S}}^{1}$ on the
measurable space $G_{\mathcal{S}}^{2}$ and $\beta$ of
$G_{\mathcal{S}}^{2}$ on the measurable space $G_{\mathcal{S}}^{1}$
can be easily calculated : take
$s=(s_{p})\in\mathcal{K}_{\mathcal{S}}$ and
$g=(g_{p})\in\mathcal{A}^{*}_{\mathcal{S}}$ such that for all
$p\in\mathcal{S}$, $g_{p}(s_{p}-1)+1\neq 0$ and, for $p$ large
enough, $g_{p}(s_{p}-1)+1\in 1+p\Zset_{p}$. Then
\begin{equation}\label{action}
\alpha_{g}(s)=\left(g_{p}(s_{p}-1)+1\right),\quad\beta_{s}(g)=\left(\frac{g_{p}s_{p}}{g_{p}(s_{p}-1)+1}\right).
\end{equation}
We define on $G_{\mathcal{S}}^{1}$ the Haar measure $\mu_{1}$
obtained, through the identification with
$\mathcal{A}^{*}_{\mathcal{S}}$, from the Haar measure
$\mu^{\times}$ on $\mathcal{A}^{*}_{\mathcal{S}}$. Also, we define
on $G_{\mathcal{S}}^{2}$ the Haar measure $\mu_{2}$  corresponding
to the product of the measures $\mu_{p}$ on $\Qset_{p}^{*}$, where $\mu_{p}$
is the Haar measure on $\Qset_{p}^{*}$ such that $\mu_{p}(1+p\Zset_{p})=1$. Taking
into account equation $(\ref{action})$, we see that $\alpha$ is a
restricted direct product action for $p\in\mathcal{S}$ of the
$\alpha^{p}\,:\,\Qset_{p}^{*}\rightarrow\text{Aut}(\Qset_{p}^{*})$,
$\alpha^{p}_{g_{p}}(s_{p})=g_{p}(s_{p}-1)+1$. Also $\beta$ is a
restricted direct product action of
$\beta^{p}\,:\,\Qset_{p}^{*}\rightarrow\text{Aut}(\Qset_{p}^{*})$,
$\beta^{p}_{s_{p}}(g_{p})=\frac{g_{p}s_{p}}{g_{p}(s_{p}-1)+1}$. We introduce the notation $\nu_{p}$ for the Haar measure on $\Qset_{p}$ such that $\nu_{p}(\Zset_{p}^{*})=1$. We have 
\begin{displaymath}
\mu_{p}=(p-1)\mu_{p}^{\times}\, ,\,\,\nu_{p}=(1-p^{-1})^{-1}\mu_{p}^{+}\,\,\text{and}\,\, \d\mu_{p}^{+}(x)=(1-p^{-1})|x|_{p}\d\mu_{p}^{\times}(x).
\end{displaymath}

The main result of this section is the following theorem which implies the description of the types of the factors $M_{\mathcal{S}}$ and $\hat{M}_{\mathcal{S}}$.

\begin{thm}
For any infinite subset $\mathcal{S}$ of $\mathcal{P}$ we have
the following isomorphisms
\begin{displaymath}
M_{\mathcal{S}}\simeq
N_{\mathcal{S}}\qquad\text{and}\qquad\hat{M}_{\mathcal{S}}\simeq
M_{\mathcal{S}}\otimes\mathcal{R},
\end{displaymath}
where $N_{\mathcal{S}}$ is the
Boca-Zaharescu factor and $\mathcal{R}$ is the hyperfinite
$\rm{II}_{1}$ factor.
\end{thm}

\begin{pf}
Let $\pi_{p}$ be the canonical inclusion of $L^{\infty}(\Qset_{p})$ in $\Qset_{p}^{*}\ltimes L^{\infty}(\Qset_{p})$. We first prove the following lemma.

\begin{lem}\label{Vaes}
Let $\mu$ be a Haar measure on $\Qset_{p}^{*}$ and $\nu$ a Haar measure on
$\Qset_{p}$. Let $K\subset\Zset_{p}^{*}$ be a subgroup of finite index with
$\mu(K)=1$ and $L$ a compact open subset of $\Zset_{p}$ such
that $KL=L$ and $\nu(L)=1$. Define 
\begin{displaymath}
e(K,L)=(\lambda(\11_{K})\otimes
1)\pi_{p}(\11_{L})\qquad\text{and}\qquad\xi(K,L)=\11_{K\times L}.
\end{displaymath}
Then $e(K,L)$ is
a projection in $\Qset_{p}^{*}\ltimes L^{\infty}(\Qset_{p})$ and
\begin{displaymath}
\left(e(K,L)\left(\Qset_{p}^{*}\ltimes
L^{\infty}(\Qset_{p})\right)e(K,L),\omega_{\xi(K,L)}\right) \simeq
\left(\mathcal{B}(l^{2}(\Nset)),\omega\right)
\end{displaymath}
where $\omega$ is the
faithful normal state on $\mathcal{B}(l^{2}(\Nset))$ with eigenvalue
list given by
\begin{displaymath}
\frac{\mu_{p}^{+}(K)}{\mu_{p}^{+}(L)}p^{-n}\quad\text{with multiplicity}\quad\left|\quotient{L\cap
p^{n}\Zset_{p}^{*}}{K}\right|,\quad n\in\Nset.
\end{displaymath}
\end{lem}

\begin{pf}
The fact that $e(K,L)$ is a projection has been proved in
Proposition \ref{crossed}. We define the following unitary
\begin{displaymath}
U\,\,:\,\,L^{2}(\Qset_{p}^{*}\times\Qset_{p},\mu\times\nu)\rightarrow
L^{2}(\Qset_{p}^{*}\times\Qset_{p}^{*},\mu\times\mu)
\end{displaymath}
\begin{displaymath}
(U\xi)(x,y)=\left(\frac{\mu_{p}^{+}(K)}{\mu_{p}^{+}(L)}|y|_{p})\right)^{\frac{1}{2}}\xi(xy^{-1},y).
\end{displaymath}
Then
\begin{displaymath}
U\lambda_{g}\otimes 1 U^{*}=\lambda_{g}\otimes
1\quad\text{and}\quad U\pi_{p}(F) U^{*}=F\otimes 1,
\end{displaymath}
this implies
that 
\begin{displaymath}
U\Qset_{p}^{*}\ltimes L^{\infty}(\Qset_{p})
U^{*}=\mathcal{B}(L^{2}(\Qset_{p}^{*},\mu))\otimes 1.
\end{displaymath}
Next, we have
\begin{displaymath}
(e(K,L)\xi)(x,y)=\11_{L}(xy)\int_{\Qset_{p}^{*}}\11_{K}(t)\xi(t^{-1}x,y)\d\mu(t)
\end{displaymath}
thus, after a simple computation, we obtain
\begin{displaymath}
Ue(K,L)U^{*}=f(K,L)\otimes 1\quad\text{where,}
\end{displaymath}
\begin{displaymath}
(f(K,L)\xi)(x)=\11_{L}(x)\int_{\Qset_{p}^{*}}\11_{K}(t)\xi(t^{-1}x)\d\mu(t).
\end{displaymath}
Observe that the image of $f(K,L)$ is the set of functions $\xi\in
L^{2}(\Qset_{p}^{*},\mu)$ such that the support of $\xi$ is in $L-\{0\}$ and
$\xi$ is invariant under translations of $K$. Writing
\begin{displaymath}
L-\{0\}=\cup_{n\in\Nset}L\cap p^{n}\Zset_{p}^{*},
\end{displaymath}
we see that every function $\xi$ in the image of $f(K,L)$ is of the
form
\begin{displaymath}
\xi=\sum_{n\in\Nset}\sum_{[y]\in\quotient{L\cap
p^{n}\Zset_{p}^{*}}{K}}\xi(y)\11_{[y]}.
\end{displaymath}
Thus we have
\begin{displaymath}
f(K,L)L^{2}(\Qset_{p}^{*},\mu)=\overline{\text{Span}<\11_{[y]},\,\,\,n\in\Nset,\,\,\,[y]\in\quotient{L\cap
p^{n}\Zset_{p}^{*}}{K}>},
\end{displaymath}
where  $\overline{\text{Span}<X>}$ means the closed vector space generated by $X$. Because $\mu([y])=\mu(K)=1$, the set of vectors $\11_{[y]}$ for $[y]\in\quotient{L\cap
p^{n}\Zset_{p}^{*}}{K}$ and $n\geq 0$ is an orthonormal basis of $f(K,L)L^{2}(\Qset_{p}^{*},\mu)$. Thus, there is a unitary $W$ between $f(K,L)L^{2}(\Qset_{p}^{*},\mu)$ and $l^{2}(\Nset)$ such that 
\begin{displaymath}
(W\otimes 1)Ue(K,L)\left(\Qset_{p}^{*}\ltimes L^{\infty}(\Qset_{p})\right)e(K,L)
U^{*}(W^{*}\otimes 1)=\mathcal{B}(l^{2}(\Nset))\otimes 1,
\end{displaymath}
and, using the computation 
\begin{displaymath}
U\xi(K,L)=\sum_{n}\sum_{[y]\in\quotient{L\cap
p^{n}\Zset_{p}^{*}}{K}}\lambda_{n,[y]}^{\frac{1}{2}}\11_{[y]}\otimes\11_{[y]},
\end{displaymath}
where $\lambda_{n,[y]}=\frac{\mu_{p}^{+}(K)}{\mu_{p}^{+}(L)}p^{-n}$, we conclude the proof.\qed
\end{pf}

\begin{rem}
We obtain, for $(K,L,\mu,\nu)=(\Zset_{p}^{*},\Zset_{p},\mu_{p}^{\times},\mu_{p}^{+})$, the list $(1-p^{-1})p^{-n}$ with multiplicity one and, for $(K,L,\mu,\nu)=(1+p\Zset_{p},\Zset_{p}^{*} -1,\mu_{p},\nu_{p})$, the following list : $(p-1)^{-1}$ with multiplicity $p-2$ and $(p-1)^{-1}p^{-n}$ with multiplicity $p-1$ for $n\geq 1$.
\end{rem} 

The next ingredient of the proof is the following lemma.

\begin{lem}\label{Vaes2}
For any infinite subset $\mathcal{S}\subset\mathcal{P}$ we have
\begin{enumerate}
\item
$M_{\mathcal{S}}\simeq\bigotimes_{p\in\mathcal{S}}\left(\Qset_{p}^{*}\ltimes L^{\infty}(\Qset_{p}),L^{2}(\Qset_{p}^{*}\times\Qset_{p},\mu_{p}^{\times}\times\mu_{p}^{+}),\xi(\Zset_{p}^{*},\Zset_{p})\right),$
\item
$\hat{M}_{\mathcal{S}}\simeq\bigotimes_{p\in\mathcal{S}}\left(\Qset_{p}^{*}\ltimes L^{\infty}(\Qset_{p}),L^{2}(\Qset_{p}^{*}\times\Qset_{p},\mu_{p}\times\nu_{p}),\xi(1+p\Zset_{p},\Zset_{p}^{*}-1)\right).$
\end{enumerate}
\end{lem}

\begin{pf}
To obtain the first isomorphism, recall that 
\begin{displaymath}
G_{\mathcal{S}}^{1}\ltimes L^{\infty}(G_{\mathcal{S}}^{2})\simeq G_{\mathcal{S}}^{1}\ltimes L^{\infty}(\quotient{G_{\mathcal{S}}}{G_{\mathcal{S}}^{1}}),
\end{displaymath}
and because $G_{\mathcal{S}}^{1}=\mathcal{A}^{*}_{\mathcal{S}}\times\{0\}$, it is easy to see that 
\begin{displaymath}
G_{\mathcal{S}}^{1}\ltimes L^{\infty}(G_{\mathcal{S}}^{2})\simeq\mathcal{A}^{*}_{\mathcal{S}}\ltimes L^{\infty}(\mathcal{A}_{\mathcal{S}}).
\end{displaymath}
Next, using Proposition \ref{crossed}, we obtain immediately the first isomorphism. For the second isomorphism, we first use Proposition \ref{crossed} and the discussion preceding the lemma to obtain 
\begin{displaymath}
\hat{M}_{\mathcal{S}}\simeq
\bigotimes_{p\in\mathcal{S}}\left(\Qset_{p}^{*}\,_{\beta^{p}}\ltimes L^{\infty}(\Qset_{p}^{*}),L^{2}(\Qset_{p}^{*}\times\Qset_{p}^{*},\mu_{p}
\times\mu_{p}^{\times}),\11_{(1+p\Zset_{p})\times\Zset_{p}^{*}}\right).
\end{displaymath}
Now define
\begin{displaymath}
V\,:\,L^{2}(\Qset_{p}^{*}\times\Qset_{p}^{*},\mu_{p}\times\mu_{p}^{\times})\rightarrow
L^{2}(\Qset_{p}^{*}\times\Qset_{p},\mu_{p}\times\nu_{p}),
\end{displaymath}
\begin{displaymath}
(V\xi)(g,s)=\xi(g^{-1},(s+1)^{-1})|s+1|_{p}^{\frac{1}{2}}.
\end{displaymath}
$V$ is unitary and
\begin{displaymath}
V\left(\Qset_{p}^{*}\, _{\beta^{p}}\ltimes L^{\infty}(\Qset_{p}^{*})\right) V^{*}=\Qset_{p}^{*}\ltimes L^{\infty}(\Qset_{p}),
\end{displaymath}
where the action for the crossed product on the right is the translation. Finally, the computation
\begin{displaymath}
V\11_{(1+p\Zset_{p})\times\Zset_{p}^{*}}=\11_{(1+p\Zset_{p})\times(\Zset_{p}^{*}-1)}
\end{displaymath}
concludes the proof.\qed
\end{pf}

We can now prove the Theorem. Using Lemmas \ref{Vaes} and \ref{Vaes2}, the remark between these two lemmas and Proposition \ref{qMq} we obtain, using the notation $e_{1}=\bigotimes_{p\in\mathcal{S}}e(\Zset_{p}^{*},\Zset_{p})$ and $e_{2}=\bigotimes_{p\in\mathcal{S}}e(1+p\Zset_{p},\Zset_{p}^{*} -1)$,  
\begin{displaymath}
e_{1}M_{\mathcal{S}}e_{1}\simeq N_{S}\quad\text{and},\quad
e_{2}\hat{M}_{\mathcal{S}}e_{2}\simeq 
\bigotimes_{p\in\mathcal{S}}\left(\mathcal{B}(l^{2}(\Nset)),\psi_{p}\right),
\end{displaymath}
where the eigenvalue list of $\psi_{p}$ is given by 
\begin{displaymath}
(p-1)^{-1}p^{-n}\quad\text{with multiplicity}\quad
\left\{ \begin{array}{ll}
p-2 & \text{if $n=0$} \\
p-1 & \text{if $n\geq 1$.}
\end{array} \right.
\end{displaymath}
Next, because $M_{\mathcal{S}}$, $\hat{M}_{\mathcal{S}}$, $e_{1}$ and $e_{2}$ are purely infinite and $e_{1}$ and $e_{2}$ have central support equal to $1$ we have that $e_{1}M_{\mathcal{S}}e_{1}\simeq M_{\mathcal{S}}$ and $e_{2}\hat{M}_{\mathcal{S}}e_{2}\simeq\hat{M}_{\mathcal{S}}$. Thus, to conclude the proof, it is sufficient to prove that 
\begin{displaymath}
\bigotimes_{p\in\mathcal{S}}\left(\mathcal{B}(l^{2}(\Nset)),\psi_{p}\right)\simeq N_{\mathcal{S}}\otimes\mathcal{R}.
\end{displaymath}
Using Lemma \ref{valprop} and  
\begin{displaymath}
\sum_{p\in\mathcal{S}}\sum_{n\geq 1}(p-1)^{-1}p^{-n}=\sum_{p\in\mathcal{S}}(p-1)^{-2}<\infty
\end{displaymath}
we can remove one copy of $p^{-1},p^{-2},\ldots$ without changing the isomorphism class of the ITPFI factor (the projection obtained in Lemma \ref{valprop} is clearly purely infinite), thus, we obtain 
\begin{displaymath}
\bigotimes_{p\in\mathcal{S}}\left(\mathcal{B}(l^{2}(\Nset)),\psi_{p}\right)\simeq\bigotimes_{p\in\mathcal{S}}\left(\mathcal{B}(l^{2}(\Nset))\otimes M_{p-2}(\Cset),\omega_{p}\otimes\tau_{p}\right),
\end{displaymath}
where  $\tau_{p}$ is the normalized trace of the matrix algebra $M_{p-2}(\Cset)$.  The theorem follows.\qed
\end{pf}

\begin{cor}
For any infinite subset $\mathcal{S}\subset\mathcal{P}$, we have
\begin{enumerate}
\item $\sum_{p\in\mathcal{S}}\frac{1}{p}<+\infty \Leftrightarrow \mu^{+}(\mathcal{A}_{\mathcal{S}}-\mathcal{A}^{*}_{\mathcal{S}})=0 \Leftrightarrow (M_{\mathcal{S}},\Delta_{\mathcal{S}})$ is of type $(\rm{I}_{\infty},\rm{II}_{\infty})$.\label{cor1}
\item $\sum_{p\in\mathcal{S}}\frac{1}{p}=+\infty \Leftrightarrow \mu^{+}(\mathcal{A}^{*}_{\mathcal{S}})=0 \Leftrightarrow (M_{\mathcal{S}},\Delta_{\mathcal{S}})$ is of type $(\rm{III},\rm{III})$.
\end{enumerate}
Moreover we have
\begin{itemize}
\item For any $\lambda\in [0,1]$ there exists a subset $\mathcal{S}\subset\mathcal{P}$ such that $(M_{\mathcal{S}},\Delta_{\mathcal{S}})$ is of type $(\rm{III}_{\lambda},\rm{III}_{\lambda})$.
\item For any countable subgroup $K$ of $\Rset$ and countable subset $\Sigma$ of $\Rset -K$ there exists a subset $\mathcal{S}$ of $\mathcal{P}$ such that $T(M_{\mathcal{S}})$ contains $K$ and does not intersect $\Sigma$.
\end{itemize}
\end{cor}

\begin{pf}
We have $\mu^{+}\left(\prod_{p\in\mathcal{S}}\Zset_{p}^{*}\right)=0 \Leftrightarrow
\mu^{+}(\mathcal{A}^{*}_{\mathcal{S}})=0$. Then, because
\begin{displaymath}
\mu^{+}\left(\prod_{p\in\mathcal{S}}\Zset_{p}^{*}\right)
=\prod_{p\in\mathcal{S}}\mu_{p}^{+}\left(\Zset_{p}^{*}\right)
=\prod_{p\in\mathcal{S}}\left(1-\frac{1}{p}\right),
\end{displaymath}
we have
\begin{equation}\label{equiv}
\sum_{p\in\mathcal{S}}\frac{1}{p}=+\infty \Leftrightarrow
\mu^{+}(\mathcal{A}^{*}_{\mathcal{S}})=0.
\end{equation}
Now, the Borel-Cantelli lemma gives 
$\sum_{p\in\mathcal{S}}\frac{1}{p}<+\infty \Rightarrow
\mu^{+}(\mathcal{A}_{\mathcal{S}}-\mathcal{A}^{*}_{\mathcal{S}})=0,$
then $(\ref{equiv})$ implies that the last implication is an
equivalence. Note that for any \lc{} ring $\mathcal{A}$, such that
$\mathcal{A}-\mathcal{A}^{*}$ has additive Haar measure  zero, the
translation action of $\mathcal{A}^{*}$ on $\mathcal{A}$ is free and
ergodic, and the corresponding crossed product is a type
$\rm{I}_{\infty}$ factor, the proof of $(\ref{cor1})$ follows. 

Now, suppose that
$\sum_{p\in\mathcal{S}}\frac{1}{p}=\infty$ then
\begin{displaymath}
\sum_{p\in\mathcal{S},i\geq
0}p^{-i}(1-p^{-1})\inf\{|p^{i}-1|^{2},1\}=
\sum_{p\in\mathcal{S},i\geq
1}p^{-i}(1-p^{-1})=\sum_{p\in\mathcal{S}}\frac{1}{p}=+\infty.
\end{displaymath}
This implies, taking into account preliminaries about ITPFI factors, that $M_{\mathcal{S}}$ and $\hat{M}_{\mathcal{S}}$ are type
$\rm{III}$ factors. The last results follow from \cite{Boca}.\qed
\end{pf}

There is a minor modification of the preceding example. Take 
\begin{displaymath}
G_{\mathcal{S}}=\mathcal{K}_{\mathcal{S}}\ltimes\mathcal{A}_{\mathcal{S}}
\end{displaymath}
and define the following subgroups
\begin{displaymath}
G^{1}_{\mathcal{S}}=\mathcal{K}_{\mathcal{S}}\times\{0\}\quad\text{and}\quad G^{2}_{\mathcal{S}}=\left\{\left(a_{p},\frac{1-a_{p}}{p}\right),\,(a_{p})\in\mathcal{K}_{\mathcal{S}}\right\}.
\end{displaymath}
Then it is easy to see that $(G^{1}_{\mathcal{S}},G^{2}_{\mathcal{S}})$ is a matched pair. A direct computation gives, for $(a_{p}),\,(b_{p})\in\mathcal{K}_{\mathcal{S}}$,
\begin{displaymath}
\alpha_{(a_{p},0)}\left(b_{p},\frac{1-b_{p}}{p}\right)
=\left(a_{p}(b_{p}-1)+1,\frac{a_{p}(1-b_{p})}{p}\right)\quad
\text{and,}
\end{displaymath}
\begin{displaymath}
\beta_{\left(b_{p},\frac{1-b_{p}}{p}\right)}(a_{p},0)
=\left(\frac{a_{p}b_{p}}{a_{p}(b_{p}-1)+1},0\right).
\end{displaymath}
We can construct the bicrossed product \lc{} quantum group $(L_{\mathcal{S}},\Delta_{\mathcal{S}})$ having the following property.

\begin{prop}
For any infinite subset $\mathcal{S}$ of $\mathcal{P}$, the \lc{} quantum group $(L_{\mathcal{S}},\Delta_{\mathcal{S}})$ is self-dual and
\begin{displaymath}
L_{\mathcal{S}}\simeq N_{\mathcal{S}}\otimes \mathcal{R}.
\end{displaymath}
\end{prop}

\begin{pf}
Define the isomorphism $u\,:\,G^{1}_{\mathcal{S}}\rightarrow G^{1}_{\mathcal{S}}$ by 
\begin{displaymath}
u(a_{p},0)=\left(a_{p}^{-1},\frac{1-a_{p}^{-1}}{p}\right),
\end{displaymath}
one verifies that 
\begin{displaymath}
u\left(\beta_{\left(b_{p},\frac{1-b_{p}}{p}\right)}(a_{p},0)\right)
=\alpha_{u^{-1}\left(b_{p},\frac{1-b_{p}}{p}\right)}\left(u(a_{p},0)\right).
\end{displaymath}
Hence, interchanging $\alpha$ and $\beta$, we get an isomorphic matched pair and so an isomorphic \lc{} quantum group. To obtain the isomorphism, recall that $G_{\mathcal{S}}^{1}\ltimes L^{\infty}(G_{\mathcal{S}}^{2})\simeq G_{\mathcal{S}}^{1}\ltimes L^{\infty}(\quotient{G_{\mathcal{S}}}{G_{\mathcal{S}}^{1}})$, and because $G_{\mathcal{S}}^{1}=\mathcal{K}_{\mathcal{S}}\times\{0\}$, it is easy to see that $G_{\mathcal{S}}^{1}\ltimes L^{\infty}(G_{\mathcal{S}}^{2})\simeq\mathcal{K}_{\mathcal{S}}\ltimes L^{\infty}(\mathcal{A}_{\mathcal{S}})$. Next, using Lemma \ref{crossed} we obtain
\begin{displaymath}
L_{\mathcal{S}}\simeq\bigotimes_{p\in\mathcal{S}}\left(\Qset_{p}^{*}\ltimes L^{\infty}(\Qset_{p}),L^{2}(\Qset_{p}^{*}\times\Qset_{p},\mu_{p}\times\mu_{p}^{+}),\xi(1+p\Zset_{p},\Zset_{p})\right).
\end{displaymath}
This implies, using Lemma \ref{Vaes} with $(K,L,\mu,\nu)=(1+p\Zset_{p},\Zset_{p},\mu_{p},\mu_{p}^{+})$ and Proposition \ref{qMq}, the following isomorphism
\begin{displaymath}
L_{\mathcal{S}}\simeq\bigotimes_{p\in\mathcal{S}}\left(\mathcal{B}(l^{2}(\Nset))\otimes M_{p-1}(\Cset),\omega_{p}\otimes\tau_{p}\right),
\end{displaymath}
where  $\tau_{p}$ is the normalized trace of the matrix algebra $M_{p-1}(\Cset)$.\qed
\end{pf}

\begin{ack}
I would like to thank Prof. S. Vaes for suggesting the subject of the paper and for many comments and Prof. F. Boca for helpful discussions. I am also grateful to Prof. L. Vainerman for many helpful discussions and comments.
\end{ack}

\bibliographystyle{plain}
\bibliography{PFima.bib}

\end{document}